\documentclass[12pt]{article}

\begin{document}

\author{Isaac M. Sonin \\
Dept. of Mathematics, UNC Charlotte\\
Charlotte, NC 28223\bigskip}
\title{Independent Events in a Simple Random Experiment and the Meaning of Independence }
\date{\today}
\maketitle

\begin{abstract} 
We count the number and patterns of pairs and tuples of independent events in a simple random experiment: first a fair coin is flipped and then a fair die is tossed. The first number, equal to 888,888, suggest that there are some open questions about the structure of independence even in a finite sample space. We discuss briefly these questions and possible approaches to answer them.

\end{abstract}




\textbf{Introduction.} \smallskip In \cite
{kolm33} A. Kolmogorov wrote: ``The concept of mutual \textit{independence} of two or more experiments holds, in a certain sense, a central position in the theory of probability.'' However, although the definition of independence is straightforward, it seems that even for a very simple
sample space it is not quite clear what the ``sources'' of independence are and how to describe them.

Let us consider the following simple random experiment: first we flip a coin and then we
toss a die. Our sample space consists of $12$ outcomes each having a
probability of $1/12$. This experiment is used in many textbook as an
illustration of the concept of independent events. However, as far as we know,
no one has asked the following two simple questions.

\textbf{Question 1.} How many different pairs $(A,B)$ of independent events
are there ?

\textbf{Question 2.} How many different tuples of independent events $%
(A_{1},A_{2},...,A_{k})$ are there ?

Let the numbers $K_{1}$ and $K_{2}$ be the answers to these questions. They can be calculated exactly, and are a bit puzzling. The first number $K_{1}$ is equal to a rather strange looking number $888,888,$
which may seem too large. On the other hand, the number $K_{1}$ may also seem too
small. It is only about 10\% of the number of all possible unordered
pairs $(A,B)$ of subsets of this sample space, which is $\frac{1}{2}\ast
2^{12}\ast (2^{12}+1)=8,390,656.$ The second number, $K_{2}=30,826,488.$ As
we will see, most of these pairs and tuples are isomorphic and can be
obtained in a small number of different ways or patterns.
This raises one more question: are the answers to Question 1 and Question 2 counterintuitive or not ? and
what do they suggest ? It seem they suggest two open problems.

\textbf{Problem 1.} How many different \textit{patterns} of independence are there
for any (finite) sample space and how to describe them ?

 In fact, when we think about
independent events, we have in mind two different concepts covered by the
same definition. The first is: independent events are those which are
produced by two or more independent sources of randomness, i.e. independent
``random generators.'' The second is: information that B has occurred
does not changed the probability of A.
So it seems that there is another open question.

\textbf{Problem 2. }How to distinguish the two concepts mentioned above?

In the next section we calculate the numbers $K_{1}$ and $K_{2}$, and in the
last section we briefly discuss Problems 1 and 2. We are far
from giving conclusive answers to these questions and hope that this small
note will encourage other probabilists to express their opinions. We also note that there is a huge literature treating these questions but we did not find direct relations to our example and even a brief survey will substantially increase the size of this note.
\newline
The author would like to thank Harold Reiter, Ernst Presman and Michael Grabchak who read the
first version of this paper and made valuable comments.

\textbf{Answers to Questions 1 and 2. }Let us denote $|A|=a,|B|=b,$ and $%
|AB|=d,$ where $AB=A\cap B.$ Since the order in the pair $(A,B)$ does not
matter, without loss of generality we can assume that $a\leq b.$ Formally,
two events are always independent if one of them is the whole sample space
or an empty set but this is not ''real independence'' so we will assume that
$0<a,b<12.$ Then it is easy to see that $1\leq d<a,b$. While potentially $d=1,2,...,10$, we will see later that $d\leq 6,$ (see also Proposition 1
in section 3).

If $d=1,$ then $\frac{1}{12}=\frac{a}{12}\frac{b}{12},$ i.e. $ab=12=3\ast
4=2\ast 6.$ If $a=3,b=4$ then there are
${12 \choose 1}$
choices for $AB,$
${11 \choose 2}$ for the two remaining elements of $A,$ and
${9 \choose 3}$
choices for the three remaining elements of $B.$ Similar reasoning works
for the case when $a=2,b=6.$ Hence the number of such pairs is
${12 \choose 1}{11 \choose 2}{9 \choose 3}=55,440=n_{1},$ and
${12 \choose 1}{11 \choose 1}{10 \choose 5}=12\ast 11\ast 252=33,264=n_{2}.$

If two events $A$ and $B$ are independent then the pairs $A$ and $B^{C}
$, $\ A^{C}$ and $B$, and $A^{C}$ and $B^{C}$ are also
independent. This means that if $a=3,b=4$ then we will have three more cases.
Slightly abusing notation, and keeping the assumption that $a\leq b,$ we can
describe them as: 2) $a=3,b=8,$ $d=2;$ 3)  $a=4,b=9,$ $d=3;$ and 4) $a=8,b=9,
$ $d=6$. In all four cases the number of pairs is the same,  $n_{1}=55,440$. We can describe all of them by an (ordered)
partition $N_{1}=[1,2,3,6].$ The multinomial coefficient $n_{1}$ is the total number of such partitions. Each element of any partition specifies a unique pair of independent events.

For the case $a=2,b=6,d=1$ we have only one more case $a=6,b=10,d=5,$
because in this case $|B|=|B^{C}|=6.$ We can describe these two cases
by a partition $N_{2}=[1,1,5,5]$. The number of pairs in each case is the multinomial coefficient $n_{2}.$

If $d=2,$ then $\frac{2}{12}=\frac{a}{12}\frac{b}{12},$ i.e. $ab=24=3\ast
8=4\ast 6.$ We do not consider other factorizations since they will not
satisfy other restrictions, e.g. $24=2\ast 12$ since $d=2<a,b$. Hence the
number of such pairs is $n_{1}$ in the first case and
$n_{3}={12 \choose 2}{10 \choose 2}{8 \choose 4}=66\ast 45\ast 70=207,900$ in the second case.
The number $n_{3}$ corresponds to a partition $N_{3}=[2,2,4,4]$, which produces
one more case $a=6,b=8,$ $d=4.$

If $d=3,$ then $\frac{3}{12}=\frac{a}{12}\frac{b}{12},$ i.e. $ab=36=4\ast
9=6\ast 6.$ Again, we must skip some other factorizations. The factorization
$36=4\ast 9$ corresponds to a partition $N_{1}$. The factorization $36=6\ast
6$ corresponds to a partition $N_{4}=[3,3,3,3]$, which produces the number
of pairs $n_{4}=\frac{1}{2}{12 \choose 3}{9 \choose 3}{6 \choose 3}=\frac{%
1}{2}(220\ast 84\ast 20)=184,800.$ Note that the factor $\frac{1}{2}$ is
present because in a pair of sets $(A,B)$ both have the same size 6.

If $d=4,$ then $ab=48=6\ast 8.$ The other factorizations
are impossible. This is a case from partition $N_{3}$. Hence
the number of such pairs is $n_{3}.$ \medskip

If $d=5,$ then similarly $ab=5\ast 12=60=6\ast 10.$ This is a case from partition $%
N_{2}$. Hence the number of such pairs is $n_{2}.$

\medskip If $d=6,$ then similarly $ab=6\ast 12=72=8\ast 9.$ This is a case from partition $N_{1}$. Hence the number of
such pairs is $n_{1}.$ \medskip

If $d\geq 7,$ then $ab=12d$ and it is easy to see that there no possible
factorizations since $d<a,b<12.$

In summary we have the following table of all possible pairs

\medskip
\begin{tabular}{|l|l|l|l|l|l|l|}
\hline
$d;ab$ & 1 ; 12 & 2 ; 24 & 3 ; 36 & 4 ; 48 & 5 ; 60 & 6 ; 72 \\ \hline
$d\ast b=$ & $3\ast 4,2\ast 6$ & $3\ast 8,4\ast 6$ & $4\ast 9,6\ast 6$ & $%
6\ast 8$ & $6\ast 10$ & $8\ast 9$ \\ \hline
$Partition$ & $N_{1}$\hspace{0.1in} ;\hspace{0.1in} $N_{2}$ & \hspace{0.1in}
$N_{1}$ ;\hspace{0.1in} $N_{3}$ & \hspace{0.1in}$N_{1}$ ;\hspace{0.1in} $%
N_{4}$ & \hspace{0.1in} $N_{3}$ & \hspace{0.15in} $N_{2}$ & \hspace{0.1in} $%
N_{1}$ \\ \hline
\end{tabular}
\newline

And the total is number $K_{1}=4n_{1}+2n_{2}+2n_{3}+n_{4}=888,888.$

\textbf{Answer to Question 2. } If $|ABC|=e=1,$ and $|A|=a,|B|=b,|C|=c$\
then $\frac{1}{12}=\frac{a}{12}\frac{b}{12}\frac{c}{12},$ $abc=144=3\ast
3\ast 2^{4}=6\ast 6\ast 4.$ Since pairs $(A,B),(A,C)$ and $(B,C)$ must also
be independent, we consider only the decompositions of $abc$ into
factors such that corresponding pairs are present in our Table. Thus, though e.g.
$144=3\ast 6\ast 8=9\ast 4\ast 4,$ since there is no pairs $3\ast
6$, or $4\ast 4$ in our table, we do not analyze them. A triplet $A,B,C$
with $e=1,$ and $a=b=6,c=4,$ implies automatically through $\frac{6}{12}\ast
\frac{6}{12}=\frac{3}{12}$ and other similar equalities that $%
|AB|=3,|AC|=|BC|=2,$ and a partition $N_{5}=[1,1,1,1,2,2,2,2]$. Then there are ${12 \choose 1}$ choices for $ABC,$
${11 \choose 5}$ for the five remaining elements of $A,$ ${5 \choose 2}$
choices for the two elements of $B$ which are in $ABC^{C},{6 \choose 3}$
choices for the three elements of $B$ outside of $A,$ ${3 \choose 1}$ choices
for the element of $C$ in $AB^{C}C,$ ${3 \choose 1}$ choices for the element
of $C$ in $A^{C}BC,$ and ${3 \choose 1}$ choices for the element of $C$
outside of $A\cup B.$ We must also divide the product by two because sets $A$
and $B$ both have the same size 6 and their order does not matter. The total
number of such triplets is \medskip $n_{5}={12 \choose 1}{11 \choose 5}{5 \choose 2}{6 \choose 3}{3 \choose 1}{3 \choose 1}{3 \choose 1}\ast
\frac{1}{2}=12\ast 462\ast 10\ast 20\ast 3\ast 3\ast 3\ast \frac{1}{2}%
=14,968,800.$

If $e=2,$\ then $abc=2\ast 144=3\ast 3\ast 2^{5}=6\ast 6\ast 8=8\ast 9\ast
2=...$ Products $6\ast 6$ and $6\ast 8$ are present in our table but $2\ast 9
$ and other are not so we need only to analyze any triplet with $a=b=6,c=8,$
and hence automatically with $|AB|=3,|AC|=|BC|=4$. This is the other case
for the partition $N_{5}$ above, and hence the total number of such triplets is again
\medskip $n_{5}={12 \choose 2}{10 \choose 4}{4 \choose 1}{6 \choose 3}{3 \choose 2}{3 \choose 2}{3 \choose 2}\ast \frac{1}{2}=66\ast 210\ast 4\ast
20\ast 3\ast 3\ast 3\ast \frac{1}{2}=14,968,600.$

If $e=3,$\ then $abc=3\ast 144=3\ast 3\ast 3\ast 2^{4}=6\ast 6\ast 12=8\ast
9\ast 6=....$ For any factorization there are pairwise products which are not
in our table. If $e>3$ then $abc=e\ast 3\ast 3\ast 4\ast 4$ and then at least one of $a,b,c$ must be at least 12. Similar reasoning shows that there is no four or more independent events. \\
Thus the total number of independent tuples,
including pairs, is $K_{2}=K_{1}+2n_{5}=30,826,488.$

Note that if $|S|=11,$ or $13,$ and $p(s)=1/|S|$ then for such sample spaces
there are no independent events at all. The following remark is due to E.
Presman:

\textbf{Proposition 1.} If a sample space $S$ consists of $n$ equally likely
outcomes and there are $k$ independent events with $r$ common points then
\begin{equation}
{0<\frac{r}{n}\leq ([n/2]/n)^{k}.}
\end{equation}
The proof follows from the equality $P(A_{1}...A_{k})=P(A_{1})...P(A_{k})$,
where the left side is equal to $\frac{r}{n},$ and in the right side each set $%
A_{i}$ can be replaced by its complement if necessary.

\textbf{Discussion}. We have found four different patterns of independence in our example. We do not have a definitive answer to Problem 1 in the general case. A possible approach to Problem 2 is as follows. If we slightly change the probabilities
of the sample points, the independence will disappear for almost all pairs,
and it is easy to show that there are arbitrarily small perturbations
which result in the complete absence of independent events. On the other hand, if
we change arbitrarily the probabilities of the independent "generators", i.e.if, in our case, we
have a biased coin and a biased die then about one hundred events will stay
independent anyway, i.e.they are ``truly independent''. A well-known example due to S. Bernstein gives three dependent events, such that any pair of these events is independent; such independence is unstable ! The
same is true of the following simpler example. There are three coins, two of them are fair and are flipped independently, the third coins is placed H(ead)
up, if the first coin and the second show H, T or T, H, and placed T(ail)
otherwise. The three coins are not independent but each pair of coins is.

Similar questions about the structure of random generators and
independent events make sense if we consider \textit{conditionally}
independent events. They correspond to random generators triggered by the
outcomes of other random generators.

To conclude this note we give another citation from a subsection on independence in \cite
{kolm33} : ``In consequence, one of the most important problems in the
philosophy of the natural sciences is -- in addition to the well-known one
regarding the essence of the concept of probability itself -- to make precise
the premises which would make it possible to regard any given real events as
independent. This question, however, is beyond the scope of this book.''

\end{document}